\documentstyle[12pt]{article}
\begin{document}

\begin{center}
{\Large\sc On exponential sums}\\
$\ $ \\
{by Ricardo Garc\'{\i}a L\'opez \footnote{Partially supported by
the DGCYT, PB95-0274.}}
\end{center}

\newcommand{\ea}{{\bf A}^{n}}
\newcommand{\al}{{\bf A}^{1}}
\newcommand{\ti}[1]{\tilde{#1}}
\newcommand{\cpx}{{\bf C}}
\newcommand{\prom}{{\bf P}}
\newcommand{\pr}{k[x_1,\ldots,x_{n}]}
\newcommand{\hpr}{K \{ x_1,\ldots,x_{n} \} }
\newcommand{\prac}{K[x_1,\ldots,x_{n}]}
\newcommand{\prc}{\bar{{\bf k}}[x_1,\ldots,x_{n}]}
\newcommand{\pe}{{\bf P}^{n}}
\newcommand{\ppe}{{\bf P}^{n-1}}
\newcommand{\xd}{X_{f}^d}
\newcommand{\ch}[2]{H^{#2}(#1)}
\newcommand{\chhc}[3]{H^{#1}_c(#2,#3)}
\newcommand{\chh}[3]{H^{#1}(#2,#3)}
\newcommand{\chc}[2]{H^{#1}_c(#2)}
\newcommand{\xf}{X^{\infty}_f}
\newcommand{\cp}{{\bf F}_p}
\newcommand{\sh}[1]{{\cal F}(#1)}
\newcommand{\gf}{f^{-1}(\bar \eta )}
\newcommand{\cgf}{{\bar f}^{-1}(\bar \eta )}
\newcommand{\va}[1]{ \{ #1=0 \} }
\newcommand{\mg}[1]{H^{#1}_c ({\bf A}^n, {\cal F}(f))}
\newcommand{\natm}{{\bf N}}
\newcommand{\ql}{{\bf Q}_l}
\newcommand{\vs}{\vspace{6mm}}
\renewcommand{\labelenumi}{\roman{enumi})}

\noindent {\bf 0.\  Introduction.}

\vs

\noindent (0.1) Let $k$ denote a finite field with $q=p^s$
elements, let $f\in\pr$ be a polynomial and let
$\Psi:\cp \to \cpx ^{\ast}$ be a non-trivial additive character.
Consider the exponential sum:
\[
S(\Psi, f)=\sum_{x\in k^n}\Psi(\mbox{Tr}_{k/\cp}(f(x))).
\]
This sum admits the following cohomological interpretation (see
\cite[Sommes trig.]{SGA4.5}, \cite[(3.5.4)]{ES}):
Let $H_f\subseteq {\bf
A}_k^{n+1}=\mbox{Spec}k[x_1,\dots,x_n,t]$ be the hypersurface given by
the
equation $t^p-t=f$. The projection map $H_f\to \ea _k$ gives a Galois
covering of $\ea _k$ (an Artin-Schreier covering) with Galois group the
additive group $\cp$.
Fix a prime $l\neq p$ and a finite extension $E_{\lambda}$ of
${\bf Q}_l$ containing the $p$-roots of unity, fix an inmersion ${\bf
Q}(\mu _p)\hookrightarrow E_{\lambda}$. One can regard $\Psi$ as a
character on $\cp $ with values on $E^{\ast}_{\lambda}$ and, extending
the structure group of the torsor $H_f$ by means of $\Psi $, one obtains
an $E_{\lambda}$-sheaf on $\ea _k$ which will be denoted ${\cal
F}_k(f)$. From now on, a scheme over $k$ (or a sheaf defined on such a
scheme) will be denoted by a subindex $k$, the supression of this
subindex will mean that we have extended scalars to a fixed algebraic
closure of $k$ that will be denoted $K$, thus we have
$\sh{f}:={\cal F}_k(f)\otimes_{k}K$. Notice that if we take the
polynomial $x\in k[x]$ and we set ${\cal L}={\cal F}(x)$, then for any
polynomial $f$ (regarded as a map $f:\ea \to \al$) one has that ${\cal
F}(f)=f^{\ast}{\cal L}$.

\noindent The cohomological interpretation of $S(\Psi, f)$ follows
now from Grothendieck's trace formula (cf.
\cite[pg. 174]{SGA4.5}). One has:
\[
S(\Psi, f)= \sum_{i}(-1)^i \mbox{Tr}(\mbox{F}, \chc{i}{\ea, \sh{f}}),
\]
where $H_c^{\ast}$ denotes cohomology with compact supports
and Tr denotes the trace of the Frobenius
morphism $\mbox{F}:\chc{i}{\ea, \sh{f}}\to\chc{i}{\ea, \sh{f}}$.
\vs

\noindent (0.2) Set $d=deg(f)$, let $f=f_{d}+f_{d-1}+\dots$ be the
decomposition of $f$ into homogeneous components and let $X_f^i\subseteq
\ppe $ be the projective hypersurface defined by the form $f_i$ ($1\le
i\le d$). In \cite[Th\'eor\`eme (8.4)]{DeW1} (see also
\cite[(3.7)]{DeW2}), P. Deligne applies
his solution of the Weil conjecture to prove the following theorem:
\vs

\noindent {\bf Theorem} (Deligne): {\it With the notations above,
assume: \begin{enumerate}
\item $X_f^d$ is non-singular.
\item gcd(d,p)=1.
\end{enumerate}
\noindent Then:
\begin{enumerate}
\item $\mg{i}=0$ \ if \ $i\neq n$.
\item $\mbox{dim}\mg{n}=(d-1)^n$ \ and $\mg{n}$ is pure of weight $n$
(i.e., all eigenvalues of the Frobenius action on this vector space have
absolute value $q^{n/2}$).
\item $\| S(\Psi, f)\| \le (d-1)^n q^{n/2}$
\end{enumerate}}
\noindent Notice that iii) follows from i) and ii) in view of
Grothendieck's trace formula. According to Katz (\cite[pg.151]{ES}),
this theorem answers a question posed by Mordell and later also by
Bombieri. \vs

\noindent (0.3) In this paper we give bounds for the exponential sums
$S(\Psi ,f)$ in cases where the variety $X_f^d$ is singular, at the
price of additional (but explicit) restrictions on the characteristic of
$k$. Before we state our main result, we first describe the
singularities that we are going to allow on $X^d_f$.
\vs

\noindent {\it Definition:} Let $g\in \pr$ be a polynomial, $\delta$ a
positive integer.
We will say that $g$ is weighted homogeneous of total degree $\delta$ if
there are
positive integers $\alpha _1,\dots, \alpha_n$ (called the weights) with
$gcd(\alpha _1, \dots, \alpha _n)=1$ and such that
for all $x_1,\dots ,x_n,\lambda \in K$
one has:
\[
g(\lambda^{\alpha_1}x_1,\dots,\lambda^{\alpha_n}x_n)=\lambda^{\delta}
g(x_1,\dots,x_n).
\]
Let $X$ be an algebraic variety over $K$, let $x\in X$ be a
closed point. We will say that $X$ has a weighted homogeneous
hypersurface singularity at $x$ (whhs for short) if there is an
isomorphism of $K$-algebras
\[
{\cal O}^h_{X,x}\simeq\frac{K \{ x_1,\dots,x_n\} }{(g)}
\]
where ${\cal O}^h_{X,x}$ denotes the henselianization of ${\cal
O}_{X,x}$ and $g$ is a weighted homogeneous polynomial.
If $\delta $ is the total degree of $g$
we will say that $\delta$ is a total
degree of the (\'etale) germ $(X,x)$. If the singularity of $X$ at $x$
is isolated (i.e. $X$ is smooth in a punctured \'etale neighborhood of
$x$), as usual we will denote by $\mu$ the Milnor number of $(X,x)$,
that is,
\[
\mu = \mbox{dim}_{K}\  \frac{K \{ x_1,\dots,x_n\} }{(\partial
g /\partial x_1, \dots ,\partial g /\partial x_n)}.
\]
One has to be aware of the fact that, unlike
in the complex case,
in general the number $\mu$ is not the dimension of a space of vanishing
cycles associated to the smoothing of $(X,x)$ given by the equation $g$.

\noindent
The Milnor number depends of course on the germ $(X,x)$, but not on
the polynomial $g$. This follows from  \cite[Expos\'e
XVI, (1.1) and (1.3)]{SGA7}) and also from the following lemma,
which we recall here since it will be (implicitely) used several times
in this paper:
\vs

\noindent {\bf Lemma:} {\it Let $K$ be a field, let $g,h\in\hpr$
be such that one has an isomorphism of $K$-algebras
\[
\frac{\hpr}{(g)} \simeq \frac{\hpr}{(h)}.
\]
Then there is a $K$-automorphism $\Phi$ of $\hpr$ such that
$\Phi (g)=h$.}
\vs

\noindent The proof given in the complex case by Loojienga in
\cite[Chap. I, Lemma (1.7)]{Lo}
works verbatim for an arbitrary field $K$.

\noindent As in the case of complex singularity germs, whhs include a
number of interesting families of singularities.
For example, it follows from \cite{Ar} that if $char(K)>5$
then all rational double points are whhs (since in this case they
are given by the same normal forms as in the complex case).

\noindent In this paper we prove the following:
\newpage

\noindent (0.4)\ {\bf Theorem:} {\it Let $f\in\pr$ be a polynomial
and $\Psi:\cp \to \cpx^{\ast}$ a non-trivial additive character. With
the notations described above, assume:
\begin{enumerate}
\item $X^d_f$ has at most weighted homogeneous
isolated singularities. If $\mbox{Sing}\ (\xd)= \{ x_1,\dots ,x_s \}$,
let $\mu _i$ denote the Milnor number of $(\xd ,x_i)$ and let $\delta
_i$ be a total degree for this germ $(1\le i \le s)$.
\item $x_i\not\in X^{d-1}_f$ for $1\le i \le s$.
\item $gcd(p,d(d-1)\delta_1 \dots \delta_s)=1$.
\end{enumerate}
Then:
\begin{enumerate}
\item $\mg{i}=0$ if $i\neq n$.
\item $\mbox{dim}\mg{n}=(d-1)^n-\sum_{i=1}^{s}\mu_i$\
and $\mg{n}$ is pure of weight $n$.
\item
$
\| S(\Psi, f) \| \le ((d-1)^n-\sum_{i=1}^{s}\mu_i)\cdot q^{n/2}.
$
\end{enumerate}}
\vs

\noindent (0.5) {\it Examples:}

\noindent a) {\it Sums in two variables:}
Let $f\in k[x,y]$ an homogeneous polynomial of degre $d$. In $K[x,y]$
we have a decomposition
\[
f_d(x,y)=\prod_{i=1}^l (\alpha _i x - \beta _i y)^{n_i}.
\]
Let $f\in k[x,y]$ be any polynomial of highest degree form $f_d$
and such that if $f_{d-1}$ is its homogeneous component of degree
$d-1$ then $f(\beta _i, \alpha _i)\neq 0$ whenever $n_i>1$.
Then, if gcd$(p,d(d-1)\Pi_{i=1}^l n_i)=1$ one gets the bound
\[
\| S(\Psi, f) \| \le ((d-1)^2-\sum_{i=1}^l(n_i-1)) \cdot q.
\]

\noindent b) {\it Arrangements of lines:}
Let $l_1,\dots,\l_d\in k[x,y,z]$
be distinct linear forms defining an arrangement of lines in ${\bf P}_k
^2$.
Set $f_d=l_1 \cdot \dots \cdot l_d$ and let $f \in k[x,y,z]$ be any
polynomial with highest degree form $f_d$ and such that its homogeneous
component of degree $d-1$ defines a plane curve
that does not pass through any of the intersection points of the lines.

\noindent Let $n_i$ denote the number of points where exactly $i$ lines
meet.
Then, if gcd$(p, d(d-1)\prod_{n_i>0}i)=1$,
one gets from (0.4) the bound:
\[
\| S(\Psi, f) \| \le ((d-1)^3-\sum_{i\ge 1}n_i(i-1)) \cdot q^{3/2}.
\]
In particular, if the arrangement is generic (i.e. no more than two
lines meet at a point), one has:
\[
\| S(\Psi, f) \| \le ((d-1)^3 - \frac{(d-1)(d-2)}{2}) \cdot q^{3/2}
\]
under the assumption that gcd$(p,d(d-1))=1$.
\vs

\noindent Other extensions of Deligne's theorem are proved
in \cite[(5.1.1)]{ES} and \cite{AS}, the last generalized in
\cite[Theorem (9.2)]{DeLo}.
These theorems neither imply nor are implied by theorem (0.4) above.
The proof of (0.4) is partially inspired by the results in \cite{GaNe}
about
the monodromy at infinity of polynomials (with complex coefficients).
The rest of the paper is devoted to it.
\vs

\noindent {\bf 1. Reduction to the case $f=f_d+x_n^{d-1}$.}
\vs

\noindent (1.1) Let $f \in \pr$ be a
polynomial satisfying the assumptions of (0.4).
Let $k\subseteq k'$ be a finite extension of $k$ such that in $\ppe
_{k'}$ there is a hyperplane
$H\subseteq \ppe$ which intersects transversally
the projective hypersurface $X^d_f$. Let $\varphi $ be a linear
automorphism of $k'[x_1,\dots,x_n]$ which sends $H$ to the hyperplane
$x_n=0$. Let
$f'=\varphi (f)=f'_d+\dots $ and set $g=f'_d+x_n^{d-1}$.
\vs

\noindent We claim
that if $\chc{i}{\ea,{\cal F}(g)}$ is pure of weight $n$, so is
$\chc{i}{{\bf A}^n, {\cal F}(f')}$, and
these $E_{\lambda}$-vector spaces have the same dimension
for all $i\ge 0$.
In order to prove this claim,
we adapt to our situation a construction of
Deligne (cf. \cite[(8.10)]{DeW1}, \cite[(3.7.2)]{DeW2}):
Let $S'$ be the affine space
over $k'$ which parametrices polynomials in $n$ variables of degree
$d$ with fixed highest degree form $f'_d$. Let $S\subseteq S'$ be the
open subset corresponding to polynomials $h$ such that
$\mbox{Sing}(X^d_{h})\cap X^{d-1}_h =\emptyset$. Let $F_{S}\in
H^0(S,{\cal O}_{S}[x_1,\dots ,x_n])$ be the universal polynomial over
$S$. It defines a corresponding Artin-Schreier sheaf ${\cal
F}(F_{S})$ over $S\times {\bf A}_{k'}^n $, let $\pi:S\times {\bf
A}_{k'}^n\to{\bf A}_{k'}^n$ be the projection map. Then one has:
\vs

\noindent (1.2) {\bf Lemma:} {\it The sheaves $R^i\pi_{!}{\cal
F}(F_{S})$ are locally constant for all $i\ge 0$.}
\vs

\noindent {\it Proof:} It follows closely that of \cite[(3.7.3)]{DeW2},
we detail the proof as long as differences appear: We embed $\ea _{k'}
\hookrightarrow {\bf P}_{k'}^n$, let $V_{\infty}$ be the hyperplane at
infinity and denote by $j:S\times {\bf A}_{k'}^n \hookrightarrow S\times
{\bf P}_{k'}^n$ the inclusion map. We claim that, locally for the
\'etale topology, the $S$-scheme $S\times \pe _{k'}$ endowed with the
sheaf
$j_{!}{\cal F}(F_S)$ is constant, that is, isomorphic to the product of
$S$ with a scheme endowed with a sheaf. Over $S\times\ea _{k'}$ this is
clear because $j_!{\cal F}(F_S)$ is locally constant there. If $h\in
V_{\infty}$ is a closed point and $x$ a closed point on $\{ h \} \times
S$, then we have the following possibilities:
\begin{enumerate}
\item $f'_d$ does not vanish on $h$.

\noindent Then there are local coordinates $t_1,\dots $ in a
neighborhood of $x$
such that $F_S=t_1^{-d}$ (use that $gcd(d,p)=1$). Since this expression
for $F_S$ is independent of the parameters of $S$, the claim follows.

\item $f'_d$ vanishes at $h$, but $h$ is not a singular point of
$f'_d=0$.

\noindent Then there is an \'etale S-morphism $\varphi $ from an \'etale
neighborhood of $p$ into $S\times \ea _{k'}$ such that $F_S\circ \varphi
= t_1^{-d}t_2$, where $t_1,t_2,\dots$ are local coordinates on that
\'etale neighborhood.

\item $h$ is a singular point of $f'_d=0$.

\noindent In an affine neighborhood of
$p$ we have:
\[
F_S=(\alpha _d +t\alpha_{d-1}+t^2\alpha_{d-2}+\dots)t^{-d}
\]
where $t,\dots $ are local coordinates, $\alpha_{d}\in m_x^2,\
\alpha_{d-1}\not\in m_x$ \ ($m_x$ being the
maximal ideal corresponding to $x$) and $\alpha _d$ is independent of
the parameters of $S$ (since the highest degree form is fixed in all
polynomials corresponding to points of $S$).
Let $f_1=\alpha _d t^{-d}$, $f_2=t^{-d+1}(\alpha _{d-1}+t\alpha _{d-2}+
\dots )$.
Then, in a neighborhood of $x$ we have
$F_S=f_1+f_2$. From this equality follows (cf. e.g. \cite[(2.1.3)]{Lau})
that locally
$j_!{\cal F}(F_S)=j_!{\cal F}(f_1) \otimes j_!{\cal F}(f_2)$. Now it
suffices to prove that $j_!{\cal F}(f_2)$ is constant as an $S$-sheaf.
But this follows from the fact that after a change of coordinates,
one has $f_2 =t_1^{-d+1}$, $t_1$ a local coordinate.
\end{enumerate}

\noindent Now one concludes exactly as in \cite[(3.7.3)]{DeW2} by using
\cite[Th. finitude, (2.16)]{SGA4.5} and \cite[Appendix]{SGA4.5}. $\Box$
\vs

\noindent Both the polynomials $f'=\varphi(f)$ and $g$ correspond to
points
on the open set $S$. But if $s\in S$ and we denote by $F_S(s)$ the
corresponding polynomial, the fiber of $R^i\pi_!{\cal F}(F_S)$ over $s$
is precisely $\chhc{i}{\ea_{k'}}{{\cal F}(F_S(s))}$. Thus we have
that $\chhc{i}{\ea}{{\cal F}(f')}$ and $\chhc{i}{\ea}{{\cal F}(g)}$,
have the same dimension and moreover, by \cite[(1.8.12)]{DeW2}, if one
of them is pure of weight $n$ so is the other. Now from the
definition of $f'$ it follows that in order to prove our main result it
suffices to prove: \vs

\noindent (1.3) {\bf Proposition:} {\it Let $f\in\prac$ be a polynomial
of the form $f=f_d+x_n^{d-1}$ which satisfies the assumptions of
(0.4). Assume that the projective hypersurface
$\va{f_d}\subseteq\ppe$ intersects the hyperplane $x_n=0$ transversally.
Then conclusions i) and ii) of (0.4) hold for the polynomial $f$.}

\noindent This proposition will be proved in the next two sections.
\vs

\noindent {\bf 2. Ramification at infinity.}
\vs

\noindent (2.1) We refer to \cite{Ray} and \cite{ES} for the definition
of tame and wild ramification of a sheaf at a point. For future
reference we recall that a pro-p-group (i.e., a projective limit
of finite p-groups) has no quotients of order prime to $p$ and that
any continous $l$-adic representation of a pro-p-group $P$ factors
through a finite quotient of $P$ (cf. \cite[p. 515]{ST}).
\vs

\noindent Let $Z\subseteq \pe \times \al$ be the hypersurface given by
the
equation $f_d(x_1,\dots,x_n)+x_n^{d-1}x_0=sx_0^d$ where $x_0,\dots ,x_n$
are homogeneous coordinates in $\pe$ and $s$ is a coordinate in $\al$.
The projection map $Z\to \al$ is a compactification of $f:\ea \to \al$,
therefore we will denote it by $\bar f$. Then we have:
\vs

\noindent (2.2) {\bf Proposition:} {\it The sheaves $R^i\bar{f}
_{\ast}{\bf Q}_l$ are tamely ramified at infinity for all $i\ge 0$.}
\vs

\noindent {\it Proof:} Let $X\subseteq \pe \times \al$ be the
hypersurface given
by $t(f_d(x_1,\dots ,x_n)+x_{n+1}^{d-1}x_0) = x_0^d$, $\pi:X \to \al$
the projection map. It will be enough to see that the sheaves $R^i\pi
_{\ast}{\bf Q}_l$ are tamely ramified at $0\in \al$. Let $\delta :\al
\to \al$ be given by $\delta (t)=t^d$.
Denote by $X'$ be the normalization of the fiber product $X\times
_{\delta} \al$ and by $\pi ':X'\to \al$ the projection map. Since
$gcd(d,p)=1$ and the wild part of the inertia group of $0\in \al$
is a pro-p-group, by (2.1) it will be enough to show that the sheaves
$R^i\pi '_{\ast}{\bf Q}_l$ are tamely ramified at the origin.
\vspace{1mm}

\noindent Let $S=\mbox{Spec}\ {\hat {\cal O}}_{\al,0}$ be the spectrum
of the completion
of ${\cal O}_{\al,0}$, let $\mu\in{\hat {\cal O}}_{\al,0}$ be a
uniformizing parameter, so we have $S\simeq\mbox{Spec}K[[ \mu ]]$.
Denote again by $\pi ':X'\to S$ the map obtained from $\pi '$ above by
base change. Let $s$ (resp. $\eta$) denote the closed (resp. the
generic) point of $S$. Let $\bar \eta$ be a geometric point localized at
$\eta$.
Let $Y_s$ denote the fiber of $\pi '$ over $s$, $Y_{\bar
\eta}$ the fiber over $\bar \eta$.
From the Leray spectral sequence for $Y_{\bar \eta}\to X'$ one gets a
vanishing cycles spectral sequence:
\[
E^{p,q}_2 = \chhc{p}{Y_s}{R^q\Phi_{\eta}(\ql)} \Rightarrow
\chhc{p+q}{Y_{\bar \eta}}{\ql},
\]
equivariant with repect to the action of the inertia group
$I=\mbox{Gal}(\bar \eta / \eta)$. We want to show that the action of the
wild part $P\subset I$ is trivial on $\chhc{\ast}{Y_{\bar \eta}}{\ql}$.
Following Katz (cf. \cite[pp 176-180]{ES}),
in order to prove this we will show that the action of $P$ is trivial on
the sheaves
$R^q\Phi_{\eta}(\ql)$. This implies (via the spectral sequence above)
that
the action of $P$ on $\chhc{p+q}{Y_{\bar \eta}}{\ql}$ is nilpotent. Since
this action factors through a finite quotient of $P$, it must be
semisimple and therefore trivial.
\vs

\noindent Let $z\in Y_s$ be a closed point. If $z$ is not a singular
point of $Y_s$, then
\[
(R^q\Phi_{\eta}(\ql))_z \simeq
\left\{ \begin{array}{rl}
0   & \mbox{ if }  q>0 \\
\ql & \mbox{ if }  q=0
\end{array} \right.
\]
and the action of $I$ is trivial. In general one has
\[
(R^q\Phi_{\eta}(\ql))_z \simeq \chhc{q}{X'_{(z){\bar \eta}}}{\ql}
\]
where $X'_{(z){\bar \eta}}$ is the fiber over $\bar \eta$ of the
henselianization $X'_{(z)}$ of $Y$ at $z$ (the ``local Milnor fiber").
Notice that $X'$ is isomorphic to the hypersurface in $\pe _S$ given by
$\va{f_d(x_1,\dots ,x_{n+1})+\mu x_{n+1}^{d-1}x_0 -x^d_0}$.
Under this isomorphism, a (closed) singular point $z$ of $Y_s$
corresponds to a point $((0:x_1:\dots :x_n),s) \in \pe \times S$, where
$x=(x_1:\dots :x_n)$ is a
singular point of $X^d_f$. Let $g$ be a local equation
defining the germ $(X^d_f,x)$ which
is a weighted homogeneous polynomial with weights $(\alpha_1, \dots ,
\alpha_{n})$ and total degree $\delta$. Then one has
an isomorphism
\[
X'_{(z)} \simeq \mbox{Spec}\left( \frac{\hpr}{g+\mu x_0 - x_0^d}
\right)
\]
Let $T$ be an indeterminate, let $X_{j}$ \ ($j=1,2$) be the
$K(T)[[ \mu ]]$-schema obtained from $X'$ by the base change $\theta
_j:K[[\mu ]] \to K(T)[[\mu ]]$ defined by
\[
\theta _j = \left\{
\begin{array}{ll}
\mu & \mbox{ if } \ j=1 \\
T^{(d-1)\delta}\mu & \mbox{ if } \ j=2
\end{array} \right.
\]
The coordinates change given by:
\[
\left.
\begin{array}{lll}
T   & \mapsto & T \\
\mu & \mapsto & \mu \\
x_0 & \mapsto & x_0 T^{\delta} \\
x_j & \mapsto & x_j T^{d\alpha _j} \ \ j=1,\dots ,n
\end{array} \right\}
\]
establishes a $K(T)[[\mu ]]$-isomorphism between $X_{1(z_1)}$ and
$X_{2(z_2)}$, where $z_1$ (resp. $z_2$) is the only $K(T)$-rational
point of $X_1$ (resp. of $X_2$) over $z$. By \cite[(9.3.5.2)]{Lau}
(see also \cite[pp. 184-187]{ES}),
the action of the inertia group $I$ over $\chhc{\ast}{X'_{(z){\bar
\eta}}}{\ql}$ factors through the finite cyclic quotient of $I$ of order
$(d-1)\delta$. Since we assumed that gcd$((d-1)\delta ,p)=1$ and
$P$ is a pro-p-group, by (2.1) the action of $P$ is trivial and we are
done. $\Box$ \vs

\noindent (2.3) {\bf Corollary:} {\it The sheaves $R^if_{!}\ql$ are
tamely ramified at infinity for all $i\ge 0$.}
\vs

\noindent {\it Proof:}
Let $X^{\infty}=X'\cap\va{x_0}$, let $X^0=X'-X^{\infty}$, set
$\pi^0=\pi '_{\mid X^0}$, $\pi ^{\infty}=\pi '_{\mid X^{\infty}}$ (as
above, $\pi ':X'\to S$).
In order to prove the corollary it suffices to prove that the sheaves
$R^i(\pi^0)_{!}\ql$ are tamely ramified at the closed point of $S$.
There is an exact sequence of sheaves on $S$:
\[
\dots \to R^{i-1}\pi ^{\infty}_{\ast}\ql
\to R^i\pi^0_{!}\ql \to
R^i\pi'_{\ast}\ql \to  \dots
\]
equivariant with respect to the action of the inertia group. One has
$X^{\infty}\simeq \xd \times S$ and under this identification $\pi
^{\infty}$ corresponds to the projection onto $S$. This implies that the
action of the inertia group on $R^{i-1}\pi ^{\infty}_{\ast}\ql$ is
trivial.
By the proof of the proposition above, the action of the wild part
$P\subset I$ is trivial on $R^i\pi'_{\ast}\ql$. Again since
$P$ is a pro-p-group, its action on $R^i\pi ^0_{!}\ql$ is
semisimple and then, in view of the exact sequence above, is trivial.
$\Box$
\newpage

\noindent {\bf 3. Cohomology computations and end of the proof.}
\vs

\noindent We first recall the following proposition, which is a direct
consequence of Proposition (3.1) in \cite{DeLo}:
\vs

\noindent (3.1) {\bf Proposition:}\ {\it Let $Y$ be a smooth $K$-scheme
of pure
dimension $n$, $g:Z\to \al$ a proper $K$-morphism. Suppose that $g$
is smooth outside a finite number of points and $R^ig_{\ast}{\bf
Q}_l$ has tame ramification at infinity for all $i\ge 0$. Then
\begin{eqnarray*}
& \chhc{i}{Z}{g^{\ast}{\cal L}}& =  0 \ \ \mbox{for $i\neq n$,
and}\\
& \chhc{n}{Z}{g^{\ast}{\cal L}}&  \ \ \mbox{is pure of
weight $n$}
\end{eqnarray*}}

\noindent (In \cite{DeLo} it is assumed that both $Y$ and $g$ are
defined over $k$, but the proof works also over $K$).
\vs

\noindent We want to apply (3.1) to ${\bar f}:Z\to \al$ (cf.
section 2). In order to do this we prove the following two lemmas:
\vs

\noindent (3.2) {\bf Lemma:}\ {\it Let $f=f_d+x_n^{d-1}\in\prac$ be a
polynomial of degree $d$ and assume that gcd$(p,d-1)=1$
(as before, $p=\mbox{char}(K)$). If the hyperplane $\va{x_n}$ in
$\ppe$
intersects the hypersurface $\va{f_d}$ transversally, then the map
$f:\ea \to \al$ is smooth outside a finite number of points.}
\vs

\noindent {\it Proof:} Let $\Sigma \subseteq \ea$ denote the critical
variety of $f:\ea \to \al$. Let $V\subseteq \ea$ be the subscheme
defined by the equations
\[
\frac{\partial f_d}{\partial x_1}=\dots = \frac{\partial f_d}{\partial
x_{n-1}}=0.
\]
We claim that dim$V\le 1$. Since $V$ is a cone with vertex at
the origin it will be enough to see that the intersection $V\cap
\va{x_n}$ reduces to $0\in \ea$.
Let $W=\va{f_d}\cap\va{x_n}\subseteq \ppe$. Then $V\cap\va{x_n}$ is the
affine cone over the singular locus of $W$. But since the intersection
of $\va{f_d}$ and $\va{x_n}$ is transversal, $W$ is smooth and
$V\cap\va{x_n}=\{ 0 \}$.

\noindent If $V=\{ 0 \}$ we are done, otherwise it will be a union of
lines through the origin. Assume now that dim$\Sigma\ge 1$. Then the
affine hypersurface $\frac{\partial f_d}{\partial x_n}+(d-1)x_n^{d-2}=0$
contains at least one line $l\subseteq V$. This line cannot be contained
in the hyperplane $x_n=0$,
so it will be given by equations
\[ \left.
\begin{array}{rcl}
x_1     &    =   & \alpha_1 x_n \\
        & \vdots &              \\
x_{n-1} &    =   & \alpha_{n-1} x_n
\end{array} \right\},
\]
where $\alpha_1,\dots ,\alpha_{n-1} \in K$. Since $\frac{\partial
f_d}{\partial x_n}$ is homogeneous of
degree $d-1$, substituying on $\frac{\partial f}{\partial x_n}$ we have:
\[
x_n^{d-2}(x_n\frac{\partial f_d}{\partial
x_n}(\alpha_1,\dots,\alpha_{n},1)+(d-1))\equiv 0
\]
which implies $d-1\equiv 0$ (mod $p$), and this is impossible since
gcd$(p,d-1)=1$.
$\Box$
\vs

\noindent (3.3) {\it Remark:} In characteristic zero, it is easy to see
that
the conclusion of the lemma above holds under the weaker assumption that
the hyperplane $x_n=0$ does not intersect the singular locus of the
projective hypersurface $f_d=0$.

\noindent In positive characteristic, the lemma does not hold without
the assumption gcd$(p,d-1)=1$ (Consider for example the polynomial
$x^{p+1}+y^{p}\in K[x,y]$, char$(K)=p$).
\vs

\noindent (3.4) {\bf Lemma:} {\it Let $f\in\prac$ be a polynomial
veryfing
the assumptions of the preceding lemma. Let $\bar f:Z\to \al$ be the
compactification of $f$ introduced in section 2.
Then $\bar f$ has a finite number of critical points.}
\vs

\noindent {\it Proof:} Let $x=(x_0,\dots ,x_n)$ be a (closed)
critical point of $\bar f$. If $x_0\neq 0$ then $x$ corresponds to a
critical point of $f$ via
\begin{eqnarray*}
i: \ea & \hookrightarrow & Z  \\
(x_1,\dots ,x_n) & \mapsto & (1,x_1,\dots ,x_n,f(x_1,\dots ,x_n))
\end{eqnarray*}
and we know from (3.3) that there is a
finite number of them.
If $x_0=0$ then one must have $\frac{\partial f _d}{\partial
x_1}(x)=\dots
=\frac{\partial f_d}{\partial x_n}(x)=x_n=0$, which contradicts our
assumptions on $f$. $\Box$
\vs

\noindent (3.5) {\bf Lemma:}\ {\it
Let $Z^{\infty}=Z\cap\va{x_0}$. Then:
\[
\chhc{k}{Z^{\infty}}{\bar f^{\ast} {\cal L}_{\mid Z^{\infty}}}=0
\ \ \mbox{ for all }  k\ge 0.
\] }
\noindent {\it Proof:}
We have
that $Z^{\infty}\simeq
X_f^d\times \al$, where $X_f^d=\va{f_d} \subseteq \ppe$.
and under this identification
$\bar f _{\mid Z^{\infty}}$ corresponds to the projection onto $\al$.
One has:
\[
\chhc{k}{Z^{\infty}}{\bar f^{\ast} {\cal L} _{\mid Z^{\infty}}}=
\chhc{k}{Z^{\infty}}{(\bar f _{\mid Z^{\infty}})^{\ast}{\cal L}}=
\chhc{k}{X^d_f \times \al}{(\bar f_{\mid Z^{\infty}})^{\ast} {\cal L}}=,
\]
and by the K\"unneth formula,
\[
= \oplus_{i+j=k}\chhc{i}{X^d_f}{E_{\lambda}}\otimes
\chhc{j}{\al}{{\cal L}}
\]
But since the character $\Psi$ is non-trivial, $\chhc{j}{\al}{{\cal
L}}=0$ for all $j\ge 0$ (\cite[pg. 99, Lemme]{ES}, thus we are done.
$\Box$
\vs

\noindent (3.6) We have an
exact sequence:
\[
\dots \to \chhc{i}{\ea}{i^{\ast}\bar f^{\ast}{\cal L}} \to
\chhc{i}{Z}{\bar f^{\ast}{\cal L}} \to
\chhc{i}{Z^{\infty}}{\bar f ^{\ast}{\cal L}_{\mid Z^{\infty}}}\dots .
\]
Since $\bar f \circ i =f$, we have $i^{\ast}\bar f^{\ast}{\cal L}={\cal
F}(f)$, thus from proposition (3.1), lemmas (3.4), (3.5) and this
sequence it follows that $\chhc{i}{\ea}{{\cal F}(f)}=0$ if $i\neq n$ and
$\chhc{n}{\ea}{{\cal F}(f)}$ is pure of weight $n$. It remains to
compute its dimension.
If $X$ is a $K$-scheme and ${\cal F}$ is a $\ql$-sheaf over $X$,
set:
\[
\chi _c (X,{\cal F})= \sum_{i\ge 0}\mbox{dim}_{\ql}\chhc{i}{X}{{\cal
F}}. \]
We will use the following proposition:
\vs

\noindent (3.7){\bf Proposition} (\cite[pg. 156]{ES}):
{\it Let $f:\ea \to \al$ be a morphism such that $R^if_{!}\ql$ has tame
ramification at infinity for all $i\ge 0$. Then:
\[
\chi_c(\ea , f^{\ast}{\cal L})=1-\chi _c (f^{-1}({\bar
\eta}), \ql ),
\]
where $\bar \eta$ is a generic geometric point of $\al$.}
\vs

\noindent (3.8) Assume now that $f\in\prac$ is a polynomial satisfying
the assumptions of (1.3). Notice that $\xd=Z^{\infty} \cap \cgf$.
From the exact sequence
\[
\dots \to \chhc{i}{\cgf}{\ql} \to
\chhc{i}{\xd}{\ql}\to\chhc{i+1}{\gf}{\ql}\to\dots
\]
follows that $\chi _c(\gf ,\ql) = \chi _c (\cgf , \ql) -\chi
_c(\xd,\ql)$.
\noindent Also, $\cgf$ is a smooth projective hypersurface of degree $d$
in $\pe$, thus its Euler characteristic can be computed, say by
comparison methods, and it is well known to be:
\[
\chi _c (\cgf , \ql)= \frac{1}{d}[(1-d)^{n+1}-1]+n+1.
\]
Thus it remains only to compute the Euler characteristic of
$\xd$.
\vs

\noindent (3.9) {\bf Lemma:}\ {\it In the situation described above, and
with the same notations,
\[
\chi _c (\xd, \ql) =
\frac{1}{d}[(1-d)^{n}-1]+n+(-1)^n\sum_{j=1}^s \mu_j.
\]
where the $\mu _j$ \ $(1\le j\le s)$ are the Milnor numbers (cf.
(0.3)) of the singular points of $\xd $.}

\vs

\noindent {\it Proof:} Let $\varphi:X\to S$ be a flat, proper
$K$-morphism,
where $X$ is a smooth $K$-scheme, $S$ is the spectrum of a
strictly henselian ring with residue field $K$ and generic point $\eta$,
the fiber of $\varphi
$ over a generic geometric point $\bar \eta$ is a smooth projective
hypersurface
of degree $d$ in $\ppe$ and the special fiber is isomorphic to
$\xd$ (The existence of such a $\varphi$ can be proved as
follows: One considers the universal family ${\cal U}\to{\bf P}^N$ of
smooth projective hypersurfaces of degree $d$ in $\ppe$, the
hypersurface
$\xd$ corresponds to a point $x^{\infty}$ on the discriminant
locus, one takes then a ``transversal slice",
i.e. $S=\mbox{Spec}\ ({\cal O}_{\ppe ,x^{\infty}}/I)^{sh}$
($sh$ denotes strict henselianization), where
$I\subset {\cal O}_{\ppe, x^{\infty}}$ is an ideal
defining a smooth curve in ${\bf P}^N$ which intersects transversally
the
discriminant locus at $x^{\infty}$. The map $\varphi$ is then obtained
by base change, details are left to the reader).
\vs

\noindent In this situation we have a exact sequence of vanishing cycles
(cf. \cite[Expos\'e XIII, (2.1.8.9)]{SGA7}):
\[
\dots \to \chh{i}{\xd}{\ql} \to \chh{i}{X_{\bar \eta}}{\ql}\to
{\bf H}^i(\xd, R\Phi _{\bar \eta}(\ql))\to\dots.
\]
and
\[
{\bf H}^i(\xd, R\Phi _{\bar \eta}(\ql))=\oplus_{j=1}^s
(R^i\Phi _{\bar \eta}(\ql))_{\bar x_j},
\]
where ${\bar x_1},\dots ,{\bar x_s}$ are geometric points over the
singular points of $\xd$. Thus,
\[
\chi _c(\xd ,\ql)=\chi _c(X_{\bar \eta},\ql) -
\sum_{j=1}^s\sum_{i}
(-1)^i \mbox{dim}\ (R^i\Phi _{\bar \eta}(\ql))_{{\bar x}_j}.
\]
On the spaces
$
(R^i\Phi _{\bar \eta}(\ql))_{{\bar x}_j}
$
there is an action of the inertia group $I=\mbox{Gal}(\bar \eta /\eta)$
and, since the singularities of $\xd$ are weighted homogeneous
with total degrees $\delta _1,\dots ,\delta _s$ and gcd$(p,\delta
_1\dots\delta_s)=1$,
one proves as in proposition (2.3) that this action is moderate (i.e.,
that the action of the wild part $P\subset I$ is trivial). It follows
then from \cite[Expos\'e XVI, Th\'eor\`eme 2.4]{SGA7} that
\[
(-1)^n\sum_i(-1)^i R^i\Phi_{\bar \eta}(\ql)_{{\bar x}_j}=\mu_j,
\]
where $\mu_j$ is the Milnor number of the germ $(\xd ,x_j)$. Now
the Euler characteristic of $X_{\bar \eta}$ can be easily computed as
in (3.8) above (notice that the dimension of $X_{\bar \eta}$ is one
less than that of ${\bar f}^{-1}({\bar \eta})$), and then the lemma
is proved. $\Box$ \vs

\noindent (3.10) It follows from (3.7), (3.8) and (3.9) that
\[
\mbox{dim}\ \chhc{n}{\ea}{{\cal F}(f)}=(d-1)^n-\sum_{i=1}^s\mu_i\ .
\]
In view of (1.3), this ends the proof of the main theorem (0.4).

\begin{footnotesize}

\end{footnotesize}
\vs

\noindent {\it Departamento de Algebra y Geometr\'{\i}a, Universidad de
Barcelona. Gran V\'{\i}a, 585. E-08007 Barcelona, Spain. e-mail adress:
rgarcia@cerber.mat.ub.es}


\begin{thebibliography}{10}

\bibitem{AS}
A.~Adolphson and S.~Sperber
\newblock {Exponential sums and Newton polyedra: Cohomology and
estimates}.
\newblock {\em Annals of Math.}, 130:367--406, 1989.

\bibitem{Ar}
M.~Artin.
\newblock {Coverings of the Rational Double Points in Characteristic p}.
\newblock In W.L. Baily and T.~Shioda, editors, {\em {Complex Analysis and
  Algebraic Geometry}}, pages 11--22. Cambridge, 1977.

\bibitem{DeW1}
P.~Deligne.
\newblock {La conjecture de Weil. I}.
\newblock {\em Publ. Math. IHES}, 43:273--307, 1974.

\bibitem{DeW2}
P.~Deligne.
\newblock {La conjecture de Weil. II}.
\newblock {\em Publ. Math. IHES}, 52:137--252, 1980.

\bibitem{SGA7}
P.~Deligne and N.~Katz.
\newblock {\em Groupes de Monodromie en G\'eometrie Alg\'ebrique, SGA 7 II},
  volume 340 of {\em Lecture Notes in Mathematics}.
\newblock Springer Verlag, 1973.

\bibitem{DeLo}
J.~Denef and F.~Loeser.
\newblock {Weights of Exponential sums, Intersection Cohomology and Newton
  Polyedra}.
\newblock {\em Invent. Math.}, 106:275--294, 1991.

\bibitem{SGA4.5}
P.~Deligne et~al.
\newblock {\em Cohomologie Etale, SGA 4 1/2}, volume 569 of {\em Lecture Notes
  in Mathematics}.
\newblock Springer Verlag, 1977.

\bibitem{GaNe}
R.~Garc\'{\i}a and A.~N\'emethi.
\newblock {On the Monodromy at Infinity of a Polynomial Map}.
\newblock {\em Compos. Math.}, 100:205--231, 1996.

\bibitem{ES}
N.~Katz.
\newblock {\em Sommes exponentielles}, volume~79 of {\em Ast\'erisque}.
\newblock Soc. Math. de France, 1981.

\bibitem{Lau}
G.~Laumon.
\newblock Majoration de sommes exponentielles attach\'ees aux hypersurfaces
  diagonales.
\newblock {\em Ann. scient. Ec. Norm. Sup., 4e. serie}, 16:1--58, 1983.

\bibitem{Lo}
E.~Loojienga.
\newblock {\em Isolated Singular Points on Complete Intersections}, volume~77
  of {\em London. Math. Soc. Lecture Notes Series}.
\newblock Cambridge University Press, 1984.

\bibitem{Ray}
M.~Raynaud.
\newblock {Caract\'eristique d'Euler-Poincar\'e d'un faisceau et cohomologie
  des vari\'et\'es ab\'eliennes}.
\newblock S\'eminaire Bourbaki 1964/65, expos\'e 286.

\bibitem{ST}
J.P.~Serre and J.~Tate.
\newblock {Good reduction of abelian varieties}.
\newblock {\em Annals of Math.}, 88:492--517, 1968.

\end{thebibliography}
\end{document}